\documentclass[11pt,reqno]{amsart}
\usepackage{latexsym}

\newtheorem{theorem}[equation]{Theorem}
\newtheorem{lemma}[equation]{Lemma}

\numberwithin{equation}{section}

\theoremstyle{definition}
\newtheorem{definition}[equation]{Definition}

\newtheorem{remark}[equation]{Remark}
\newtheorem*{remark*}{Remark}

\newcommand{\Z}{{\mathbb Z}}

\newcommand{\frg}{{\mathfrak g}}
\newcommand{\frh}{{\mathfrak h}}
\newcommand{\frs}{{\mathfrak s}}
\newcommand{\frt}{{\mathfrak t}}
\newcommand{\frz}{{\mathfrak z}}
\newcommand{\frj}{{\mathfrak{j}}}
\newcommand{\frgl}{{\mathfrak{gl}}}
\newcommand{\frsl}{{\mathfrak{sl}}}
\newcommand{\frst}{{\mathfrak{st}}}
\newcommand{\ov}{\overline}
\newcommand{\uJ}{{\underline{J}}}
\newcommand{\uL}{{\underline{L}}}
\newcommand{\ufrj}{{\underline{\frj}}}

 \DeclareMathOperator{\eespan}{span}
 \DeclareMathOperator{\End}{End}

 \DeclareMathOperator{\ad}{ad}
 \DeclareMathOperator{\diag}{diag}
 
 \DeclareMathOperator{\TKK}{TKK}
 \DeclareMathOperator{\JCK}{JCK}
 \DeclareMathOperator{\CK}{CK}
 \DeclareMathOperator{\uce}{uce}

\newcommand{\psl}{\mathfrak{psl}}
\newcommand{\espan}[1]{\eespan_{\mathbb F}\left\{#1\right\}}
\newcommand{\subo}{{_{\bar 0}}}
\newcommand{\subu}{{_{\bar 1}}}
\begin{document}

\title{A$(n,n)$-graded Lie Superalgebras}
\author[G. Benkart]{Georgia Benkart$^{\star}$}
\address{Department of Mathematics, University of Wisconsin,Madison, WI 53706, U.S.A.}
\email{benkart@math.wisc.edu}
\thanks{$^{\star}$ Support from NSF Grant
\#{}DMS--0245082,  NSA Grant MDA904-03-1-0068,  and the
hospitality of Universidad de Zaragoza (Ayudas para Estancias de
Investigadores de Excelencia, Convocatoria 2002) are gratefully
acknowledged.}
\author[A. Elduque]{Alberto Elduque$^{\ast}$}
\address{Departamento de Matem\'aticas, Universidad de Zaragoza,
50009 Zaragoza, Spain} \email{elduque@unizar.es}
\thanks{$^{\ast}$ Supported by the
Spanish Ministerio de Ciencia y Tecnolog{\'\i}a and FEDER (BFM
2001-3239-C03-03) and by the Diputaci\'on General de Arag\'on
(Grupo de Investigaci\'on de \'Algebra).}
\author[C. Mart{\'\i}nez]{Consuelo Mart{\'\i}nez$^{\dag}$}
\address{Departamento de Matem\'aticas, Universidad de Oviedo,
C/Calvo Sotelo, s/n, 33007 Oviedo, Spain}
\email{chelo@pinon.ccu.uniovi.es}
\thanks{$^{\dag}$ Supported by BFM 2001-3239-C03-01 and  PR 01-GE15.}

\date{\today}

\begin{abstract}
 We determine the Lie superalgebras over fields of
 characteristic zero  that are graded by the root system A$(n,n)$
 of the special linear  Lie superalgebra $\mathfrak{psl}(n+1,n+1)$.
\end{abstract}

\maketitle


\section{Introduction}  Root space decompositions and gradings  by finite root
systems
have been a fundamental tool in the study of Lie algebras since the
classification
of the finite-dimensional complex simple Lie algebras by Killing and Cartan over
a hundred years ago.
Many nonsimple, finite- and infinite-dimensional  Lie algebras over arbitrary
fields of characteristic zero exhibit
a grading by one or more finite root systems.
Important examples include the affine Kac-Moody Lie algebras, the toroidal Lie
algebras,
and the intersection matrix Lie algebras introduced by Slodowy,  to name just a
few.
Any finite-dimensional simple Lie algebra over a field $\mathbb F$ of
characteristic
zero which has an ad-nilpotent element (or equivalently by the Jacobson-Morosov
theorem,
has  a copy of $\mathfrak {sl}_2(\mathbb F)$) is graded by a finite root system.

To make this concept more precise, Berman and Moody \cite{BM}
singled out a class of Lie algebras graded by finite root systems.
This class\footnote{together with the class of Lie algebras
graded by the nonreduced root systems BC$_r$ in \cite{ABG2} and
\cite{BS}} contains all of the above-mentioned algebras. The
definition given by Berman and Moody starts with a finite-dimensional split
simple
Lie algebra $\mathfrak g$ over a field $\mathbb F$ of
characteristic zero  having a root space decomposition $\mathfrak
g = \mathfrak h \oplus \bigoplus_{\alpha \in \Delta} \mathfrak
g_\alpha$ relative to a split Cartan subalgebra $\mathfrak h$.
Thus, $\mathfrak g$ is the analogue of a complex simple Lie
algebra, just as $\mathfrak{sl}_n(\mathbb F)$ is the analogue of
$\mathfrak{sl}_n(\mathbb C)$ over the field $\mathbb F$. The root
system $\Delta$ is one of the finite (reduced) root systems A$_n$,
B$_n$, C$_n$, D$_n$, E$_6$, E$_7$, E$_8$, F$_4$, or G$_2$.
Following Berman and Moody we say

\medskip
\begin{definition} \label {def1} A Lie algebra $L$ over $\mathbb F$ is graded by
the
(reduced) root system
$\Delta$ if
\begin{itemize}
\item [\rm {(1)}] $L$ contains $\mathfrak g$ as a subalgebra;
\item[\rm {(2)}]  $L = \bigoplus_{\alpha \in \Delta \cup\{0\}}L_\alpha$, where
$L_\alpha = \{x \in L \mid [h,x] = \alpha(h)x$ for all $h \in \mathfrak h\}$ for
$\alpha \in \Delta \cup \{0\}$; and \item [\rm {(3)}]  $L_0 = \sum_{\alpha \in
\Delta}[L_\alpha,L_{-\alpha}].$
\end{itemize}
\end{definition}

\medskip

Assumption (1) may be replaced by the weaker requirement that the
derivation algebra of $L$ contain a copy of $\mathfrak g$ without
making a substantial difference. If only conditions (1) and (2)
are assumed, then the ideal $L' := \left(\bigoplus_{\alpha \in
\Delta} L_\alpha \right) \oplus \left(\sum_{\alpha \in
\Delta}[L_\alpha,L_{-\alpha}]\right)$ is graded by $\Delta$, and $L$
acts as derivations on $L'$. Assumption (2) is the real heart of
the matter, for it implies that $[L_\alpha, L_\beta] \subseteq
L_{\alpha+\beta}$ (where $L_{\alpha+\beta} = 0$ if $\alpha+\beta
\not \in \Delta\cup \{0\}$), and so $L$ is graded by the abelian
group generated by $\Delta$.

The Lie algebras in Definition \ref{def1} have been classified in
(\cite{BM}, \cite{BZ}, \cite{N})  under the assumption that the center of $L$ is
trivial.
Their central extensions have been described in \cite{ABG} and their derivations
in \cite{B}.     The results in these papers have been used in an essential way
to determine the structure of the extended affine Lie algebras (see \cite{BGK},
\cite{BGKN}, \cite{AABGP}, \cite{AG}) and of the intersection
matrix Lie algebras (see \cite{BM}, \cite {BZ}).

Root decompositions also play a crucial role in  the
classification of the finite-dimensional complex simple Lie
superalgebras   (see \cite{K}).   Many nonsimple finite- and
infinite-dimensional Lie superalgebras exhibit such root
decompositions, and to better understand their structure, we
introduce a concept analogous to the one above.
A finite-dimensional simple complex Lie superalgebra  is classical
 or it is of Cartan type.   Here  we suppose that
$\mathfrak g$ is a finite-dimensional split simple classical Lie
superalgebra over a field $\mathbb F$ of characteristic zero (the
analogue of one of the classical superalgebras over $\mathbb C$).
Thus, $\mathfrak g$  has a root space decomposition $\mathfrak g =
\mathfrak h \oplus \bigoplus_{\alpha \in \Delta} \mathfrak
g_\alpha$ relative to a split Cartan subsuperalgebra $\mathfrak h$, and
the root system $\Delta$ is one of the following:

\smallskip
 \begin{itemize}
\item [\rm {(i)}]  A$({m,n})$,
\item [\rm {(ii)}] B$({m,n})$,  \hskip 3.7 truein ($\ast$)
\item [\rm {(iii)}] C$(n)$, \, D$(m,n)$, \, D$(2,1;\mu) \,
(\mu \in \mathbb F \setminus \{0,-1\})$,   F$(4)$, G$(3)$,
\item [\rm {(iv)}] P$(n)$, Q$(n).$
\end{itemize}

\smallskip
Such a Lie superalgebra $\mathfrak g$ can have a nonsplit central extension,
as  Weyl's theorem on complete reducibility fails to hold in the superalgebra
setting.
It is advantageous to take into account this behavior in formulating
the notion of a $\Delta$-graded Lie superalgebra.

\medskip

\begin{definition} \label{def2}  Let $\mathfrak g = \mathfrak h
\oplus \bigoplus_{\alpha \in \Delta} \mathfrak g_\alpha$ be a
finite-dimensional split simple classical  Lie superalgebra over a
field $\mathbb F$ of characteristic zero  relative to a split
Cartan subalgebra $\mathfrak h$ (so $\Delta$ is one of the
root systems in ($\ast$)).
 A Lie superalgebra $L$ over $\mathbb F$ is said to be
$\Delta$-graded if it satisfies the following three conditions:
\begin{enumerate}

\item[(1)] $L$ contains, as a subsuperalgebra, a central cover
$\widetilde\frg$ of $\frg$.  The subsuperalgebra $\widetilde\frg$ is called the
\emph{grading subsuperalgebra} of $L$.
\item[(2)] Let
$\overline \frg\subo$ be the unique subalgebra of $\widetilde\frg\subo$
which projects isomorphically onto $\frg\subo$ under the cover map
(see Remark \ref{barsubo} below), and
for a Cartan subalgebra $\frh$ of $\frg\subo$,  let $\overline {\frh}$ be
the Cartan subalgebra of $\overline \frg\subo$ which projects
isomorphically onto $\frh$. Then $L$ decomposes as
\begin{equation*}
L=\bigoplus_{\alpha\in\Delta\cup\{0\}}L_\alpha,
\end{equation*}
where $\Delta$ is the root system of $\frg$ relative to $\frh$,
and $L_\alpha = \{x \in L \mid [h,x] = \alpha(h)x$ for all $h \in
\overline \frh \}$ for $\alpha \in \Delta \cup \{0\}$; and
\item[(3)] $L_0=\sum_{\alpha\in\Delta}[L_\alpha,L_{-\alpha}]$.
\end{enumerate}
\end{definition}

\medskip

\begin{remark} \ When $\Delta$ is as in (ii) or (iii) of ($\ast$),
or when $\Delta$ is  of type A$(m,n)$ with $m \neq n$,  the only
central cover of $\frg$ is $\frg$ itself (see \cite{IK}). Thus,
Definition \ref{def2} coincides with the notion of a
$\Delta$-graded Lie superalgebra
 in \cite{BE1}-\cite{BE3} and \cite{GN} for those root systems.
 Garc{\'\i}a and Neher \cite{GN}
 determine  the Lie superalgebras
graded by  the finite root systems  A$_n$, B$_n$, C$_n$, D$_n$, E$_6$,
and E$_7$ over arbitrary commutative superrings by
exploiting a 3-grading on  the superalgebras and applying Jordan methods.
The definition of a $\Delta$-graded Lie superalgebra  presented in
{\cite[Defn. ~1.4]{MZ}}  can be seen to be equivalent to the one above.
Allowing a central cover of $\frg$ to
be the grading subsuperalgebra rather than just  using $\frg$ itself enables the
Cheng-Kac superalgebras to be
realized in a natural way as P$(3)$-graded superalgebras (see {\cite[Thm.
~1.2]{MZ}})
The Cheng-Kac superalgebras provide important examples of conformal
superalgebras
of finite type (\cite{CK} and \cite{GLS}).
\end{remark}

Lie superalgebras graded by the root systems B$(m,n)$ have been
described in \cite{BE2},  while the superalgebras graded by the
root systems listed in (iii) of ($\ast$) have been classified in
\cite{BE1}.   In  \cite{BE3}, Benkart and Elduque have determined
the Lie superalgebras graded by A$(m,n)$ for $m \neq n$. Included
in that work are  results on the  Lie superalgebras containing
$\frg = \mathfrak{psl}(n+1,n+1)$ and having a root space
decomposition of type A$(n,n)$ relative to a Cartan subalgebra of
$\frg$.  The complexity of the structure of such superalgebras was
another motivation for broadening the concept of a $\Delta$-graded
Lie superalgebra to the one given above. Mart{\'\i}nez and
Zelmanov \cite{MZ} have described the Lie superalgebras graded by
the root systems of type P$(n)$ and Q$(n)$.   Thus, all the
$\Delta$-graded Lie superalgebras, in the sense of  Definition
\ref{def2}, as well as all of their central extensions,  have been
determined except for the root systems of type A$(n,n)$.   It is
the goal of this present paper to  complete the classification by
treating the one remaining case of A$(n,n)$.

\medskip
Throughout  it will be assumed that the ground field $\mathbb F$
has characteristic zero.   Usually the field $\mathbb F$ will be
omitted from the notation unless it is needed for clarity. Thus, for example, we
will write
$\mathfrak{sl}_{n}$ rather than  $\mathfrak{sl}_{n}(\mathbb F)$
and $\mathfrak{sl}(n,n)$ instead of $\mathfrak{sl}(n,n)(\mathbb
F)$ or  $\mathfrak{sl}_{n,n}(\mathbb
F)$ as a notational convenience.

\bigskip

\section{Universal central extensions}
A  {\it central extension} of a Lie superalgebra $L$ is a
pair $(\widetilde L,\pi)$ consisting of a Lie superalgebra
$\widetilde L$ and a surjective Lie superalgebra homomorphism
$\pi: \widetilde L \rightarrow L$ (preserving the grading), whose
kernel lies in the center of $\widetilde L$.  If $\widetilde L$ is
perfect ($\widetilde L = [\widetilde L,\widetilde L]$), then
$\widetilde L$ is said to be a {\it cover} or {\it covering} of
$L$, and $\pi$ is referred to as the cover map or covering
homomorphism.  Any perfect Lie superalgebra $L$ has a unique (up
to isomorphism) covering ($\widehat L, \widehat \pi$) which is
universal,  called the {\it universal covering superalgebra} or
{\it the universal central extension}  of $L$.  Thus, any covering
of $L$ is isomorphic to a quotient of $\widehat L$ by a central
ideal.   A superalgebra is
 {\it centrally closed} if $\widehat L = L$.  When  the universal
 central extensions of two perfect Lie superalgebras $L_1$ and $L_2$
are isomorphic (or equivalently, when $L_1/Z(L_1) \cong
L_2/Z(L_2)$),
 the superalgebras $L_1$ and $L_2$ are said to be {\it centrally isogenous}.

The universal central extensions of the {\it basic} classical simple Lie
superalgebras  (that is, all the classical Lie superalgebras
except for P$(n)$ and Q$(n)$)   have been obtained in \cite{IK} under the
additional
hypothesis that  the kernel of the covering map is an even
subspace. But this condition is superfluous:

\begin{lemma}\label{uce}
Let $\frg$ be a basic classical simple Lie superalgebra.
\begin{enumerate}
\item[(i)] If $\frg$ is not of type {\text A}$(n,n)$, then $\frg$
is centrally closed.
\item[(ii)] If $\frg=\psl(n+1,n+1)$, then its universal central
extension $\widehat \frg = \uce(\frg)$ is given by
\begin{equation*}
\widehat \frg=\begin{cases}
 \frsl(n+1,n+1)&\text{for $n\geq 2$,}\\
 \text{``\,\hbox{\rm D}$(2,1;-1)$''}&\text{for $n=1$ \qquad (notation of
\cite{IK})}.
\end{cases}
\end{equation*}
\end{enumerate}
\end{lemma}

\begin{proof}
The universal central extension $\widehat \frg = \uce({\frg})$ is
a vector space direct sum of $\frg$ and a quotient of
$\frg\wedge\frg$, so $\dim\widehat \frg<\infty$. Let $\pi:\widehat
\frg\rightarrow\frg$ denote the projection.

Now $\widehat \frg$ is a $\frg$-module by means of $x.z=[\widehat
x,z]$ for any $x\in\frg$ and $z\in\widehat\frg$, where $\widehat x$ is
any element in $\widehat\frg$ which projects onto $x$. As a
$\frg$-module, $\widehat\frg$ has
 a composition
series having an adjoint factor and trivial factors. If $\frg$ is
not of type A$(n,n)$ for some $n \geq 1$, then by the complete
reducibility results of \cite{BE1}-\cite{BE3}, $\widehat
\frg=\ker\pi\oplus \frs$ for a suitable $\frg$-module $\frs$. But
$[\widehat \frg,\frs]=\frg.\frs\subseteq \frs$, so $\frs$ is an
ideal of $\widehat \frg$ and $\widehat\frg=[\widehat \frg,\widehat
\frg]\subseteq \frs$, since $\ker\pi$ is central. Therefore
$\ker\pi =0$.

When $\frg=\psl(n+1,n+1)$ for $n \geq 1$, then
$\frg\subo=\frsl_{n+1}\oplus\frsl_{n+1}$ is
semisimple so, as a $\frg\subo$-module, $\widehat \frg\subu
=\frs\subu \oplus (\ker\pi\cap \widehat \frg\subu)$ for some $\frg\subo$-
submodule
$\frs\subu$. But
$[\widehat \frg\subo,\frs\subu]=\frg\subo.\frs\subu\subseteq\frs\subu$
and $\widehat \frg\subu=[\widehat{\frg}\subo,\widehat{\frg}\subu]
=[\widehat{\frg}\subo,\frs\subu]\subseteq\frs\subu$. Therefore,
$\ker\pi\subseteq \widehat{\frg}\subo$, and the results of \cite{IK}
can be used to complete the proof.
\end{proof}

\medskip
\begin{remark} \ \label{barsubo} For every basic classical
Lie superalgebra $\frg$, there is a unique subalgebra $\overline
\frg\subo$ of $\widehat\frg\subo=\uce(\frg)\subo$ such that the
covering  homomorphism $\pi: \widehat \frg \rightarrow \frg$
restricts to an isomorphism $\overline \frg\subo\rightarrow
\frg\subo$,  since for $\psl(n+1,n+1)$, $\widehat \frg\subo$ is
the direct sum of
$\overline\frg\subo:=[\widehat\frg\subo,\widehat\frg\subo]$, which
is a direct sum of two copies of $\frsl_{n+1}$, and the center
$Z(\widehat\frg)$. In fact, the same result applies to any cover
of $\frg$. The even part  of the `strange' Lie superalgebra
D$(2,1;-1)$ is the direct sum of two copies of $\frsl_2$ and a
three-dimensional center.

The situation for the classical simple  Lie superalgebras $\frg$ of type
P$(n)$ and Q$(n)$ is completely analogous.  For both of them, the even
part is simple (isomorphic to $\frsl_{n+1}$), so the even part of
any cover $\widetilde \frg$ is the direct sum of its center and of
$[\widetilde \frg\subo,\widetilde \frg\subo]$, which is the required
unique subalgebra.
\end{remark}
 \medskip

 \begin{lemma}\label{perfectas}
Let $\frs$ and $\frg$ be two perfect Lie superalgebras such that
$\frs$ is a subsuperalgebra of $\frg$ and $\frs$ is centrally closed.
Let $\widetilde\frg$ be a cover of $\frg$ with
associated projection $\pi:\widetilde\frg\rightarrow\frg$. Then there
is a unique subsuperalgebra $\widetilde\frs$ of $\widetilde\frg$ such that the
restriction of $\pi$ to $\widetilde\frs$ is an isomorphism from
$\widetilde\frs$ onto $\frs$.
\end{lemma}
\begin{proof}
Because  $\pi: \pi^{-1}(\frs)\rightarrow \frs$ is a cover of $\frs$, \and $\frs$
is centrally closed,  $\pi^{-1}(\frs)=\frz\oplus \widetilde\frs$  for
subsuperalgebras $\widetilde\frs$ and $\frz$ with $\frz\subseteq
 \ker \pi \subseteq Z(\widetilde \frg)$.   Thus
$\pi\vert_{\widetilde\frs}:\tilde\frs\rightarrow \frs$ is an
isomorphism.   Moreover,
$\widetilde\frs=[\pi^{-1}(\frs),\pi^{-1}(\frs)]$,   which
demonstrates the uniqueness of $\widetilde\frs$.
\end{proof}
 \medskip

\begin{lemma}\label{Ltilde}
Let $L$ be a Lie superalgebra with an abelian subsuperalgebra $\frh$ such that
$L=\bigoplus_{\alpha\in\Delta\cup\{0\}}L_\alpha$, where $\Delta$ is a
finite subset of $\frh^*\setminus\{0\}$;  $L_\alpha=\{ x\in L \mid
[h,x]=\alpha(h)x\ \forall h\in \frh\}$ for any $\alpha\in \Delta$;
and $L_0=\sum_{\alpha\in \Delta}[L_\alpha,L_{-\alpha}]$.   Let
$\pi:\widehat L\rightarrow L$ be the universal central extension of
$L$. Then $\widehat L=\bigoplus_{\alpha\in\Delta\cup\{0\}}\widehat
L_\alpha$, where $\widehat L_\alpha =\{x\in\widehat L \mid
[h',x]=\alpha(\pi(h'))x\ \forall h'\in\pi^{-1}(\frh)\}$ for any
$\alpha\in \Delta$, and $\widehat L_0=\sum_{\alpha\in \Delta}[\widehat
L_\alpha,\widehat L_{-\alpha}]$. In particular $\ker \pi\subseteq
\widehat L_0$,  and the restriction of
$\pi$ to $\widehat L_\alpha$ is a linear isomorphism $\pi_\alpha:
\widehat L_\alpha\rightarrow L_\alpha$  for any $\alpha\in \Delta$.
\end{lemma}

\begin{proof}
The superalgebra $\widehat  L$ is an $\frh$-module in a natural way.
If for  $h\in\frh$,  \  $\mathfrak m(t)$ denotes the minimal polynomial of $\ad
h$ on
$L$, then  $\mathfrak m(t)$ divides $\prod_{\alpha\in \Delta\cup\{0\}}
\bigl(t-\alpha(h)\bigr)$, and  $\mathfrak m(\ad h') \widehat L\subseteq
Z(\widehat L)$ for any $h'\in \pi^{-1}(h)$.  Thus the minimal
polynomial of $\ad h'$ divides $t\mathfrak m(t)$.   This implies
$\widehat L=\widehat L_0\oplus\bigl( \bigoplus_{\alpha\in\Delta}\widehat
L_\alpha\bigr)$ where  the $\widehat L_\alpha$'s are as required,
and where $\widehat L_0 = \{x \in \widehat L \mid  (\ad h')^2 x = 0$ for
all $h' \in \pi^{-1}(\frh)\}$.    Now
$L' :=\bigl( \bigoplus_{\alpha\in\Delta}\widehat L_\alpha\bigr)
\oplus\bigl(\sum_{\alpha\in\Delta}[\widehat
L_\alpha,\widehat L_{-\alpha}]\bigr)$ is a perfect subsuperalgebra of
$\widehat  L$ with $\pi( L')=L$.    Since $L' + \ker \pi = \widehat L$ and
$\ker \pi \subseteq Z(\widehat L)$,
we have $\widehat L = [\widehat L, \widehat L] = [L',L'] = L'$,
and  the lemma follows.
\end{proof}

\bigskip

\section{Jordan supersystems}

Let $\Delta$ be a finite subset of an abelian group $\Gamma$ with
$0\not\in \Delta$.   Suppose  $\uJ =(J_\alpha, \alpha\in\Delta)$
is a family  of $\Z_2$-graded vector spaces  with even linear maps
$\psi_{\alpha,\beta}:J_\alpha\otimes J_\beta\rightarrow
J_{\alpha+\beta}$ whenever $\alpha,\beta,\alpha+\beta\in\Delta$,
and linear maps $\psi_{\alpha,-\alpha}:J_\alpha\otimes J_{-\alpha}
\rightarrow \bigoplus_{\beta\in\Delta}\End(J_\beta)$ for all
$\alpha\in \Delta$ such that  $-\alpha\in\Delta$. For any
$\alpha,\beta\in\Delta$ and $x_\alpha\in J_\alpha$,
$x_{-\alpha}\in J_{-\alpha}$, $x_\beta\in J_\beta$,
$\psi_{\alpha,-\alpha}(x_\alpha\otimes x_{-\alpha})(x_\beta)$ will
denote the action on $x_{\beta}$ of the $\End(J_\beta)$-component
 of $\psi_{\alpha,-\alpha}(x_\alpha\otimes
x_{-\alpha})$.

  Then $\uJ$ is said to be a \emph{Jordan
supersystem} if the direct sum
\begin{equation*}
K(\uJ)=\bigoplus_{\gamma\in\Delta\cup\{0\}} J_\gamma,
\end{equation*}
where
\begin{equation*}
J_0=\sum_{\substack{\alpha\in\Delta\\
\text{such that }-\alpha\in\Delta}}
\psi_{\alpha,-\alpha}(J_\alpha\otimes J_{-\alpha}),
\end{equation*}
becomes a Lie superalgebra under the natural product
\begin{equation*}
[x_\alpha,y_\beta] =\begin{cases} 0 & \hbox{\rm if} \
  \alpha+\beta \not \in \Delta \cup \{0\} \\
\psi_{\alpha,\beta}(x_\alpha \otimes y_\beta)
& \hbox{\rm if} \
 \alpha, \beta, \alpha+\beta \in \Delta \cup \{0\} \\
 x_\alpha(y_\beta)
& \hbox{\rm if} \
 \alpha = 0, \beta \in \Delta\\
 -(-1)^{|x_\alpha||y_\beta|}y_\beta(x_\alpha) &  \hbox{\rm if} \
 \alpha \in \Delta, \beta=0\\
 x_\alpha y_\beta - (-1)^{|x_\alpha||y_\beta|} y_\beta x_\alpha
&  \hbox{\rm if} \
 \alpha = \beta= 0. \end{cases} \end{equation*}
\noindent  The
even and odd parts of $K(\uJ)$ are  given by
\begin{equation*}
K(\uJ)\subo=\bigoplus_{\gamma\in \Delta\cup\{0\}}(J_\gamma)\subo
\quad \text{and} \quad
 K(\uJ)\subu=\bigoplus_{\gamma\in
\Delta\cup\{0\}}(J_\gamma)\subu.
\end{equation*}

Equivalently, $\uJ$ is a Jordan supersystem if it consists of the
nonzero homogeneous components of a $\Gamma$-graded Lie
superalgebra $L=\bigoplus_{\gamma\in\Delta\cup\{0\}}L_\gamma$ with
$[x_\alpha,y_\beta]=\psi_{\alpha,\beta}(x_\alpha,y_\beta)$ for
$\alpha,\beta,\alpha+\beta\in \Delta$, $x_\alpha\in L_\alpha$,
$y_\beta\in L_\beta$; and $[[x_\alpha,x_{-\alpha}],y_\beta]=
\psi_{\alpha,-\alpha}(x_\alpha\otimes x_{-\alpha})(y_\beta)$ for
$\alpha,-\alpha,\beta\in\Delta$.

\medskip
A homomorphism of Jordan supersystems  $\uJ = (J_\alpha, \alpha
\in \Delta)$ and $\uJ' = (J_\alpha', \alpha \in \Delta)$ is a
family of $\mathbb Z_2$-graded linear mappings $\eta = (\eta_\alpha, \alpha \in
\Delta)$, with $\eta_\alpha: J_\alpha \rightarrow J_\alpha'$, such
that the induced mapping $K(\uJ) \rightarrow K(\uJ')$ is a Lie
superalgebra homomorphism. Applying the same arguments as in
{\cite[Sec.~1.4]{BZ}}, we deduce:

\medskip

\noindent{\bf 1.}\quad
 There exists a universal such Lie superalgebra;   that
is, there exists a $\Gamma$-graded Lie superalgebra
$U_\uJ=\bigoplus_{\gamma\in\Delta\cup\{0\}} (U_\uJ)_\gamma$ and a
homomorphism $\nu: (J_\alpha, \alpha\in \Delta)\rightarrow
\left((U_\uJ)_\alpha, \alpha\in \Delta\right)$ of Jordan
supersystems such that for any $\Gamma$-graded Lie superalgebra
$L=\bigoplus_{\gamma\in\Delta\cup\{0\}}L_\gamma$ and any homomorphism of
Jordan supersystems $\eta: (J_\alpha, \alpha\in \Delta)\rightarrow
(L_\alpha,\alpha\in\Delta)$,  there is a unique homomorphism of
$\Gamma$-graded Lie superalgebras $\theta : U_\uJ \rightarrow L$ such
that $\theta \nu_\alpha=\eta_\alpha$ for any
$\alpha\in\Delta$. Moreover, $\nu_\alpha$ is a linear isomorphism
for every $\alpha\in\Delta$,  and $(U_\uJ)_0=\sum_{\alpha\in
\Delta}[(U_\uJ)_\alpha,(U_\uJ)_{-\alpha}]$.

\medskip

\noindent{\bf 2.}\quad Let
$L=\bigoplus_{\gamma\in\Delta\cup\{0\}}L_\gamma$ and $L'=
\bigoplus_{\gamma\in\Delta\cup\{0\}}L'_\gamma$ be two perfect
$\Gamma$-graded Lie superalgebras. If the Jordan supersystems
$(L_\alpha,\alpha\in\Delta)$ and $(L'_\alpha,\alpha\in\Delta)$ are
isomorphic, and if
$L_0=\sum_{\alpha\in\Delta}[L_\alpha,L_{-\alpha}]$ and
$L'_0=\sum_{\alpha\in\Delta}[L'_\alpha,L'_{-\alpha}]$, then $L$
and $L'$ are centrally isogenous. Moreover, if
$\uL=(L_\alpha,\alpha\in\Delta)$ is the corresponding Jordan
supersystem, then $U_\uL$ is a central cover of $L$
{\cite[Prop.~1.5]{BZ}}.

\medskip

\noindent{\bf 3.}\quad Let
$L=\bigoplus_{\gamma\in\Delta\cup\{0\}}L_\gamma$ be a
$\Gamma$-graded Lie superalgebra, and let
$\uL=(L_\alpha,\alpha\in\Delta)$ be the corresponding Jordan
supersystem.  If there exists an abelian subalgebra $\frh\subseteq
(L_0)\subo$ such that $\Delta\subseteq \frh^*$ and $L_\alpha=\{
x\in L \mid [h,x]=\alpha(h)x\ \forall h\in \frh\}$ for any
$\alpha\in\Delta\cup\{0\}$, and if $L_0=\sum_{\alpha\in
\Delta}[L_\alpha,L_{-\alpha}]$,  then $U_\uL$ is centrally closed
{\cite[Prop.~1.6]{BZ}}. In particular, $U_\uL$ is the universal
central cover of $L$.

\medskip

\noindent{\bf 4.}\quad Let
$L=\bigoplus_{\gamma\in\Delta\cup\{0\}}L_\gamma$ and
$L'=\bigoplus_{\gamma\in\Delta\cup\{0\}}L'_\gamma$ be two
$\Gamma$-graded Lie superalgebras,  and assume $(\eta_\alpha,
\alpha\in\Delta)$ is a family of  $\Z_2$-graded linear isomorphisms
$\eta_\alpha :L_\alpha\rightarrow L'_\alpha$ such that for any
$\alpha,\beta,\alpha+\beta\in\Delta$,
$\eta_{\alpha+\beta}([x_\alpha,y_\beta])
=[\eta_\alpha(x_\alpha),\eta_\beta(y_\beta)]$ for  $x_\alpha\in
L_\alpha$, $y_\beta\in L_\beta$.   If for any $\alpha\in
\Delta$,
\begin{equation*}
L_\alpha =\sum_{\substack{\delta,\gamma\in\Delta\\
\delta,\gamma\ne \pm\alpha\\
\delta+\gamma=\alpha}}[L_\delta,L_\gamma],
\end{equation*}
then $\eta$ is an isomorphism of Jordan supersystems
{\cite[Prop.~1.8]{BZ}}.

\bigskip

As a noteworthy example, let $\Delta=\{\pm 1\}\, (\subseteq
\mathbb{Z})$ and let $J$ be a Jordan superalgebra. Then
$\uJ=(J_1,J_{-1})$, where $J_1=J_{-1}=J$, is a Jordan supersystem
with
\begin{equation*}\label{psiTKK}
\begin{split}
\psi_{1,-1}(x\otimes y)(z)&=(xy)z+x(yz)-(-1)^{|x||y|}y(xz)\\
\psi_{1,-1}(x\otimes y)(w)&=-(xy)w+x(yw)-(-1)^{|x||y|}y(xw)
\end{split}
\end{equation*}
for any $x,z\in J_1$ and $y,w\in J_{-1}$. The associated Lie
algebra $K(\uJ)$ is the well-known (centerless) Tits-Kantor-Koecher
superalgebra of $J$ (see \cite[1.10--1.14]{BZ} and the references
therein), which we denote by $\TKK(J)$.

\smallskip

Also, if $\uJ=(J_1,J_{-1})$ is a Jordan pair, then $U_\uJ$
coincides with the Lie algebra $\widehat\TKK(J)$ (\emph{functorial
Tits-Kantor-Koecher algebra}) considered in \cite{S}, since both
algebras satisfy the same universal property \cite[Thm.
2.7]{S}.

\medskip

\section{A$(n,n)$-graded Lie superalgebras, $n\geq 2$}

Let $V = V\subo \oplus V\subu$  be a $\mathbb Z_2$-graded $\mathbb
F$-vector space with $\dim V\subo  = n+1 = V\subu$. Then $\End V =
(\End V)\subo \oplus (\End V)\subu$, where
 $X v \in V_{\overline {i}+\overline{j}}$ (subscript mod $2$)
 for $X \in (\End V)_{\overline i}$ and  $v \in V_{\overline j}$,
  and $\End V$ is a Lie superalgebra under the supercommutator
 product $[X,Y] = XY - (-1)^{|X|\, |Y|} YX$.
 The transformations $X$
 having supertrace $\mathfrak{str}(X) = \mathfrak{tr}_{V\subo}(X) -
\mathfrak{tr}_{V\subu}(X) = 0$
 form a subsuperalgebra.
  Choosing
 a basis of $V$ that respects the decomposition, we may
 identify that subsuperalgebra  with the Lie superalgebra
$\mathfrak{sl}(n+1,n+1)$.
 The identity matrix spans the center of $\frsl(n+1,n+1)$,  and
 factoring out the center gives the classical simple Lie superalgebra
$\psl(n+1,n+1)$.

We assume a basis for $V\subo$ consists of vectors  numbered $1,\dots, n+1$, and
one
for $V\subu$ by $\ov {1}, \dots, \overline{n+1}$,  and
number the rows and columns of matrices in $\frsl(n+1,n+1)$ accordingly by the
indices in
$\mathcal{I}=\{1,\ldots,n+1,$ $\ov {1},\ldots,\overline {n+1}\}$. Let $\frh$
denote the
Cartan subalgebra of $\frg=\psl(n+1,n+1)$
consisting of the diagonal matrices in $\frsl(n+1,n+1)$ modulo the center.
Thus, any element in $\frh$ is uniquely the class,
modulo the center, of a diagonal matrix
$\diag(a_1,\ldots,a_{n+1},a_{\ov {1}},\ldots,a_{\overline {n+1}})$ with
$a_1+\cdots+a_{n+1}=0=a_{\ov {1}}+\cdots+a_{\overline {n+1}}$.

Then the root system of $\psl(n+1,n+1)$ relative to $\frh$  is
\begin{equation} \label{A-rs}
\Delta=\{ \varepsilon_i-\varepsilon_j \mid i\ne j,\ i,j\in \mathcal I\},
\end{equation}
where $\varepsilon_i$ takes the class of a diagonal matrix
$\diag(a_1,\ldots,a_{n+1},a_{\ov {1}},\ldots,a_{\overline {n+1}})$ (with
$a_1+\ldots+a_{n+1}=0=a_{\ov {1}}+\ldots+a_{\overline {n+1}}$) to $a_i$. Thus,
$\varepsilon_1+\cdots+\varepsilon_{n+1}=0=\varepsilon_{\ov{1}}+\cdots
+\varepsilon_{\overline {n+1}}$.

\medskip

{\it  Throughout the remainder of the section it will always be assumed that
$n\geq 2$.}

\medskip

Let $\widehat \frg=\frsl({n+1},{n+1})$, which is the universal central extension
of $\frg=\psl({n+1},{n+1})$, and let $\pi:\widehat \frg\rightarrow \frg$ be the
natural projection. The unique subalgebra of $\widehat \frg\subo$ which
projects isomorphically onto $\frg\subo$ is
$[\widehat \frg\subo,\widehat \frg\subo]\cong \frsl_{n+1}\oplus\frsl_{n+1}$.
Identify
$\frh$ with the Cartan subalgebra of $[\widehat \frg\subo,\widehat \frg\subo]$
formed by its diagonal matrices (which are the diagonal matrices
satisfying $a_1+\cdots+a_{n+1}=0=a_{\ov{ 1}}+\cdots+a_{\overline {n+1}}$ as
above). Thus $\frh'=\mathbb F z \oplus\frh$ is the
Cartan subalgebra of $\widehat \frg$ consisting of the diagonal matrices,
where  $z$ is the $2(n+1) \times 2(n+1)$ identity matrix.
Extend $\varepsilon_i$ to $\frh'$ by imposing the condition $\varepsilon_i(z)=1$
for
any $i$.  Hence on $\frh'$, we have
$\varepsilon_1+\cdots+\varepsilon_{n+1}=\varepsilon_{\ov{1}}+\cdots+\varepsilon_
{\overline {n+1}}$, but this expression does not  equal $0$.

Let $L$ be an A$(n,n)$-graded Lie superalgebra. Then, by
definition, either $\frg$ or $\widehat \frg$ is the grading
subsuperalgebra in $L$. In either event, $L$ is a $\widehat
\frg$-module (this is obvious if $\widehat \frg\subseteq L$, and
otherwise, $L$ is a $\frg$-module, and hence a $\widehat
\frg$-module through $\pi$).   Write the action of $\widehat \frg$
on $L$ by $a.x$ for $a\in\widehat \frg$ and $x\in L$.   Thus,
$a.x=[a,x]$ if $\widehat \frg\subseteq L$ and $a.x=[\pi(a),x]$
otherwise.

The next lemma is fundamental in our work.  To avoid complicated
expressions, $e_{ij}$  ($i,j\in \mathcal I$) will denote  the matrix with a $1$
in
place $(i,j)$ and $0$'s elsewhere if we are working in $\widehat \frg$ or the
class of that matrix modulo the center in $\frg$, depending on the context.
Thus,
for $i\ne j$,  $\widehat \frg_{\varepsilon_i-\varepsilon_j}=\mathbb Fe_{ij}$ and
$\frg_{\varepsilon_i-\varepsilon_j}=\mathbb Fe_{ij}$ also.

\medskip

\begin{lemma}\label{ztrivial}  If $\widehat \frg$ is the grading subsuperalgebra
of
an A$(n,n)$-graded Lie superalgebra $L$,
then the central element $z$ acts trivially on $L$.   Equivalently, for
any $\alpha\in \Delta$,
\begin{equation}\label{eq:ztrivial}
L_\alpha=\{ x\in L \mid h'.x=\alpha(h')x\ \  \forall \  h'\in\frh'\}.
\end{equation}
\end{lemma}
\begin{proof}
For $1\leq i\ne j\leq {n+1}$, take $1\leq k\leq {n+1}$ different from $i$
and $j$ (which is possible since $n\geq 2$). Then $\varepsilon_i-
\varepsilon_j\pm
(\varepsilon_k-\varepsilon_{\ov {k}})\not\in\Delta$ and hence
$[e_{k\ov{k}},e_{\ov{k},k}].L_{\varepsilon_i- \varepsilon_j}=0$. But
$[e_{k\ov{k}},e_{\ov{k} k}]=e_{kk}+e_{\ov{ k}\,\ov{k}}\in \frh'\setminus\frh$,
with $(\varepsilon_i-\varepsilon_j)([e_{k\ov{k}},e_{\ov{k} k}])=0$, so
that \eqref{eq:ztrivial} holds for $\alpha=\varepsilon_i-\varepsilon_j$.
The same argument works for $\alpha=\varepsilon_{\ov{i}}-\varepsilon_{\ov{ j}}$.

Now without loss of generality, it is enough to check
\eqref{eq:ztrivial} for $\alpha=\varepsilon_1-\varepsilon_{\ov{1}}$.
If $n\geq 3$,  one can use that $\varepsilon_1-\varepsilon_{\ov{1}}\pm
(\varepsilon_2-\varepsilon_{\ov{2}})\not \in \Delta$ and proceed as above.

Therefore, for the rest of the proof, $n$ will be assumed to be
$2$. Here $\varepsilon_1-\varepsilon_{\ov{1}}
+\varepsilon_2-\varepsilon_{\ov{1}} \not\in\Delta$,
so $e_{2\ov{1}}.L_{\varepsilon_1-\varepsilon_{\ov{1}}}=0$,
while $\varepsilon_1-\varepsilon_{\ov{1}}
-(\varepsilon_2-\varepsilon_{\ov{1}})=\varepsilon_1-\varepsilon_2\in
\Delta$. Thus,
\begin{equation*}
\begin{split}
(e_{22}&+e_{\ov{1}\,\ov{1}}).(z.L_{\varepsilon_1-\varepsilon_{\ov{1}}})=
  [e_{2\ov{1}},e_{\ov{1} 2}].(z.L_{\varepsilon_1-\varepsilon_{\ov{1} }})\\
& = z.\left(e_{2\ov{1} }.(e_{\ov{1}  2}.
 L_{\varepsilon_1-\varepsilon_{\ov{1}}})\right) +
 z.\left(e_{\ov{1}  2}.(e_{2\ov{1} }.
  L_{\varepsilon_1-\varepsilon_{\ov{1}}})\right)\\
& = z.\left(e_{2\ov{1} }.(e_{\ov{1}  2}.
 L_{\varepsilon_1-\varepsilon_{\ov{1} }})\right)\\
&\subseteq z.(e_{2\ov{1} }.
 L_{\varepsilon_1-\varepsilon_2})=
 e_{2\ov{1} }.(z.L_{\varepsilon_1-\varepsilon_2})=0,
 \end{split}
\end{equation*}
by the calculation in the first paragraph of the proof.

Similarly, $(e_{11}+e_{\ov{2}\,\ov{2} }).(z.L_{\varepsilon_1-
\varepsilon_{\ov{1}}})=0$.  Hence,
since $z.x \in
L_{\varepsilon_1-\varepsilon_{\ov{1}}}$ too for any $x\in L_{\varepsilon_1-
\varepsilon_{\ov{1}}}$,
\begin{equation*}
\begin{split}
0&=\Big((e_{11}-e_{22})-(e_{\ov{1}\,\ov{1}}-e_{\ov{2}\,\ov{2}})\Big).(z.x)\\
&=(\varepsilon_1-\varepsilon_{\ov{1}})\Big((e_{11}-e_{22})-(e_{\ov{1}\,\ov{1}}-
e_{\ov{2}\,\ov{2}})\Big)(z.x)=2 \, z.x.
\end{split}
\end{equation*}
We conclude that $z.L_{\varepsilon_1-\varepsilon_{\ov{1}}}=0$ also,
thus proving the lemma.
\end{proof}
\medskip

\begin{lemma}\label{3graded}
With
\begin{equation*}
\begin{split}
&L(-1)= \bigoplus_{1\leq i, j\leq {n+1}}L_{\varepsilon_{\ov{i}}-
\varepsilon_j},\\
&L(0)= L_0\oplus \Bigg (\bigoplus_{1\leq i\ne j\leq
   {n+1}}L_{\varepsilon_i-\varepsilon_j}\Bigg) \oplus
   \Bigg(\bigoplus_{1\leq i\ne j\leq
   {n+1}}L_{\varepsilon_{\ov{i}}-\varepsilon_{\ov{j}}}\Bigg), \quad \hbox{\rm
and}\\
&L(1)= \bigoplus_{1\leq i, j\leq {n+1}}L_{\varepsilon_i-\varepsilon_{\ov{j}}},
\end{split}
\end{equation*}
the associated decomposition $L=L(-1)\oplus L(0)\oplus L(1)$ is a $3$-grading
of $L$.
\end{lemma}

\begin{proof}
The only problem is to show that $[L(1),L(1)]=0=[L(-1),L(-1)]$,
and this is clear for $n\geq 3$ since
$\varepsilon_i-\varepsilon_{\ov{j}}
+ \varepsilon_k-\varepsilon_{\ov{l}}\not \in \Delta$ for any $1\leq
i,j,k,l\leq n+1$.    However for $n=2$,
$\varepsilon_i-\varepsilon_{\ov{j}}+\varepsilon_k-\varepsilon_{\ov{l}}
=(\varepsilon_i+\varepsilon_k)-(\varepsilon_{\ov{j}}+\varepsilon_{\ov{l}})
=-(\varepsilon_r-\varepsilon_{\ov{s}})$ in case $i\ne k$ and $j\ne
l$, where $1\leq r,s\leq 3$ are such that $r\ne i,k$ and $s\ne
j,l$ (recall that
$\varepsilon_1+\varepsilon_2+\varepsilon_3=0=\varepsilon_{\ov{1}}
+\varepsilon_{\ov{2}}+\varepsilon_{\ov{3}}$). Without loss of
generality,  it suffices to prove that
$[L_{\varepsilon_1-\varepsilon_{\ov{1}}},L_{\varepsilon_2-
\varepsilon_{\ov{2}}}]=0$.

Set $\hat\Delta : =\{\varepsilon_i-\varepsilon_j \mid  i\ne j\in
\{1,2,3, \ov{1}\}\}$.  Then
\begin{equation*}
\hat L= \sum_{\alpha\in\hat\Delta}L_\alpha +
\sum_{\alpha\in\hat\Delta}[L_\alpha,L_{-\alpha}]
\end{equation*}
is a subsuperalgebra of $L$ which is $\frsl(3,1)$-graded.   To see this,
observe that $\frsl(3,1)$ embeds naturally in $\frsl(3,3)$ (using
the indices $\{1,2,3,\ov{1}\}$), and $\frsl(3,3)$ projects through
$\pi$ onto $\psl(3,3)$.   The composition gives an embedding of
$\frsl(3,1)$ into $\psl(3,3)$, because $\frsl(3,1)$ is simple.   Hence
$\frsl(3,1)$ is embedded in $L$,  and Lemma \ref{ztrivial}
guarantees that the root spaces work properly for the natural
Cartan subalgebra of $\frsl(3,1)$, which is contained in $\frh'$.

By the results in \cite[Sec.~3]{BE3}  on A$(m,n)$-graded
superalgebras for $m\ne n$,  in particular by the specific form of
the multiplication of these algebras \cite[(3.1)]{BE3}, it follows
that $L_{\varepsilon_1-\varepsilon_{\ov{1}}}
 =e_{3 \ov{1}}.L_{\varepsilon_1-\varepsilon_3}$.
 Then because $\varepsilon_1-\varepsilon_3+\varepsilon_2-\varepsilon_{\ov{2}}
 =-2\varepsilon_3-\varepsilon_{\ov{2}}\not\in \Delta$,  we have
 $[L_{\varepsilon_1-\varepsilon_3},L_{\varepsilon_2-\varepsilon_{\ov{2}}}]
 =0$, and thus  \begin{equation*}
\begin{split}
[L_{\varepsilon_1-\varepsilon_{\ov{1}}},
 L_{\varepsilon_2-\varepsilon_{\ov{2}}}] &=
 [e_{3 \ov{1}}.L_{\varepsilon_1-\varepsilon_3},
 L_{\varepsilon_2-\varepsilon_{\ov{2}}}]\\
&=[L_{\varepsilon_1-\varepsilon_3},
  e_{3\ov{1}}.L_{\varepsilon_2-\varepsilon_{\ov{2}}}].
\end{split}
\end{equation*}
Therefore, it is enough to prove that $e_{3\ov{1}}.L_{\varepsilon_2-
\varepsilon_{\ov{2}}}=0$.
But  now using the indices
$\{ 2, \ov{1},\ov{2},\ov{3}\}$, we get an embedding of $\frsl(1,3)$
in $\psl(3,3)$ as above, which shows that
 $L_{\varepsilon_2-\varepsilon_{\ov{2}}}=
 e_{2\ov{3}}.L_{\varepsilon_{\ov{3}}-\varepsilon_{\ov{2}}}$.    Hence
\begin{equation*}
\begin{split}
e_{3\ov{1}}.L_{\varepsilon_2-\varepsilon_{\ov{2}}}&=
  e_{3\ov{1}}.(e_{2\ov{3}}.
     L_{\varepsilon_{\ov{3}}-\varepsilon_{\ov{2}}})\\
  &\subseteq [e_{3\ov{1}},e_{2 \ov{3}}].
    L_{\varepsilon_{\ov{3}}-\varepsilon_{\ov{2}}} +
    e_{2\ov{3}}.
    (e_{3\ov{1}}.L_{\varepsilon_{\ov{3}}-\varepsilon_{\ov{2}}})\\
  &=0,
\end{split}
\end{equation*}
since $[e_{3\ov{1}},e_{2\ov{3}}]=0$ and
 $e_{3 \ov{1}}.L_{\varepsilon_{\ov{3}}-\varepsilon_{\ov{2}}}=0$,
because $\varepsilon_3+\varepsilon_{\ov{3}}-\varepsilon_{\ov{1}}
-\varepsilon_{\ov{2}}=\varepsilon_3+2\varepsilon_{\ov{3}}\not\in\Delta$.
\end{proof}

\bigskip

Now for any fixed $l\in \mathcal I=\{
1,\ldots,{n+1},\ov{1},\ldots,\overline{n+1}\}$,
there is
an embedding of
$\frsl({n+1},n)$  into both $\frsl({n+1},{n+1})$ and $\psl({n+1},{n+1})$
using the rows and columns indexed by the elements in
$\mathcal I\setminus\{l\}$. Consider
\begin{equation*}
\begin{split}
\Delta_{(l)}&=\{ \varepsilon_i-\varepsilon_j \mid  i\ne j,\  i,j\in \mathcal
I\setminus\{
            l\}\},  \qquad \hbox{\rm and}\\
L_{(l)}&=\Bigg({\displaystyle\bigoplus_{\alpha\in\Delta_{(l)}}L_\alpha
\Bigg )\oplus\Bigg(
     \sum_{\alpha\in\Delta_{(l)}}[L_\alpha,L_{-\alpha}]\Bigg)}\ .
\end{split}
\end{equation*}
Then $L_{(l)}$ contains the image of $\frsl({n+1},n)$ (under the
embedding above) and, thanks to Lemma \ref{ztrivial}, $L_{(l)}$ is
a A$({n},n-1)$-graded Lie superalgebra. By
{\cite[Thm.~3.10]{BE3}}, there is a unital associative
superalgebra $A_{(l)}$ such that $L_{(l)}$ is centrally isogenous
to $[\frgl({n+1},n)\otimes A_{(l)}, \frgl({n+1},n)\otimes
A_{(l)}]$.  In particular, we may assume that
\begin{equation*}
L_{\varepsilon_i-\varepsilon_j}=e_{ij}\otimes A_{(l)}
\end{equation*}
(a copy of $A_{(l)}$)  for $i\ne j\ne l\ne i$, and
\begin{equation}\label{producto}
[e_{ij}\otimes a,e_{jk}\otimes b]=(-1)^{|a|(|j|+|k|)}e_{ik}\otimes
ab
\end{equation}
for all distinct indices  $i,j,k\in \mathcal I\setminus\{ l\}$ and for all
homogeneous elements
$a,b\in A_{(l)}$.    Moreover,  $e_{ij}\otimes 1$ is the element $e_{ij}$
in the grading subsuperalgebra  (either $\frg$ or $\widehat \frg$) of  $L$.
Here
$|1|=\cdots=|n+1|=0$, while $|\ov{1}|=\cdots=|\overline {n+1}|=1$.

Under these conditions we have

\medskip

\begin{lemma}\label{l-inde} The unital associative superalgebra
$A_{(l)}$ is independent of $l$.
\end{lemma}

\begin{proof}
Assume $l\ne l'$ and take $i\ne j\in \mathcal I\setminus\{ l,l'\}$ (which can be
done since  $\mathcal I$ contains at least $6$ elements).   Then
$L_{\varepsilon_i-\varepsilon_j}=e_{ij}\otimes A_{(l)}=e_{ij}\otimes
A_{(l')}$,  so there is a linear ($\Z_2$-graded) bijection
$\varphi_{ll'}^{ij}:A_{(l)}\rightarrow A_{(l')}$ such that
$\varphi_{ll'}^{ij}(1)=1$.

However, for distinct values  $i,j,j'\in I\setminus\{l,l'\}$, we have
\begin{equation*}
[e_{ij}\otimes x,e_{jj'}\otimes
1]=(-1)^{|x|(|j|+|j'|)}e_{ij'}\otimes x
\end{equation*}
for any homogeneous $x\in A_{(l)}$.    The same  can be said for elements $x'\in
A_{(l')}$.    Consequently,
\begin{equation*}
\varphi_{ll'}^{ij}=\varphi_{ll'}^{ij'}.
\end{equation*}
That is, $\varphi_{ll'}^{ij}$ does not depend on $j$;  and in the
same way, it does not depend on $i$ either.    Therefore, there is
a linear $\Z_2$-graded isomorphism
$\varphi_{ll'}:A_{(l)}\rightarrow A_{(l')}$.

Moreover, when $i,j,k\in I\setminus\{l,l'\}$ are all different,  and
 $x,y$ are homogeneous elements of  $A_{(l)}$, we have
\begin{equation*}
[e_{ij}\otimes x,e_{jk}\otimes y]=(-1)^{|x|(|j|+|k|)}
e_{ik}\otimes xy,
\end{equation*}
and the same equation holds  for $x',y'  \in A_{(l')}$.     As a result,
\begin{equation*}
\varphi_{ll'}(xy)=\varphi_{ll'}(x)\varphi_{ll'}(y),
\end{equation*}
which implies $\varphi_{ll'}$ is an isomorphism that allows us to identify
$A_{(l)}$ and $A_{(l')}$ for any $l\neq l'$ in $\mathcal I$.
\end{proof}

\medskip
Setting  $A=A_{(1)}$ and making these identifications, we have
\begin{equation*}
L_{\varepsilon_i-\varepsilon_j}=e_{ij}\otimes A
\end{equation*}
for any $i\ne j$ in $\mathcal I$.   Moreover, the product is given
by
\begin{equation*}
[e_{ij}\otimes a,e_{jk}\otimes b]=(-1)^{|a|(|j|+|k|)}
e_{ik}\otimes ab
\end{equation*}
for all distinct indices $i,j,k \in \mathcal I$ and all homogeneous elements
$a,b\in
A$.

It follows from these considerations, that $(L_\alpha, \alpha \in
\Delta)$ (where $\Delta$ is of type A$(n,n), (n \geq 2)$ as in (\ref {A-rs}))
is a Jordan supersystem isomorphic to the Jordan supersystem
consisting of the root spaces of  the Lie superalgebra
$\frsl_{{n+1},{n+1}}(A):=[\frgl({n+1},{n+1})\otimes
A,\frgl({n+1},{n+1})\otimes A]$. The results on Jordan
supersystems in Section 3 immediately imply our main result in
this section:

\medskip

\begin{theorem}\label{ngeq3}
Let $n\geq 2$,  and let $L$ be a  Lie
superalgebra over a field $\mathbb F$ of characteristic zero graded by the root
system
\hbox{\rm A}$(n,n)$.  Then there is a unital associative superalgebra $A$
such that $L$ is centrally isogenous to
\begin{equation*}
\frsl_{{n+1},{n+1}}(A)=[\frgl({n+1},{n+1})\otimes A,\, \frgl({n+1},{n+1})\otimes
A].
\end{equation*}
Conversely, for any unital associative superalgebra $A$, any Lie
superalgebra centrally isogenous to $\frsl_{{n+1},{n+1}}(A)$ is
\hbox{\rm A}$(n,n)$-graded.
\end{theorem}

\medskip

For the converse, note that by Lemma \ref{Ltilde}, any central extension of
an  A$(n,n)$-graded Lie superalgebra is also A$(n,n)$-graded,
as are its central quotients.

\medskip

\begin{remark} \ The universal central extension of the Lie
superalgebra $\frsl_{{m+1},{n+1}}(A)$, for $m+n\geq 3$,  has been
determined in \cite[Thm.~2]{MP}. It is the so called Steinberg Lie
superalgebra $\frst_{{m+1},{n+1}}(A)$.
\end{remark}

\medskip

\section{A$(1,1)$-graded Lie superalgebras}

In this final section, let $\frg=\psl(2,2)$ and let $\widetilde\frg$ be any
central cover of
$\frg$. The unique subalgebra of $\widetilde\frg\subo$ which projects
isomorphically onto $\frg\subo$ is
$[\widetilde\frg\subo,\widetilde\frg\subo]\cong \frsl_2 \oplus\frsl_2$.
Adopting the same notation as in the previous section,  we have that the root
system
here is
\begin{equation*}
\Delta=\{ \pm 2\varepsilon_1, \pm 2\varepsilon_{\ov{1}}, \pm
\varepsilon_1\pm \varepsilon_{\ov{1}}\}
\end{equation*}
since $\varepsilon_1+\varepsilon_2=0=\varepsilon_{\ov{1}}+\varepsilon_{\ov{2}}$.
Moreover, the root spaces are given by
\begin{gather*}
\frg_{2\varepsilon_1} =\mathbb Fe_{12},\quad \frg_{-2\varepsilon_1}=\mathbb
Fe_{21},
 \quad \frg_{2\varepsilon_{\ov{1}}}=\mathbb Fe_{\ov{1}\,\ov{2}},\quad
 \frg_{-2\varepsilon_{\ov{1}}}=\mathbb Fe_{\ov{2}\,\ov{1}},\\
\frg_{\varepsilon_1+\varepsilon_{\ov{1}}}= \mathbb Fe_{1\ov{2}}\oplus \mathbb
Fe_{\ov{1} 2},
 \qquad \frg_{-(\varepsilon_1+\varepsilon_{\ov{1}})}
     =\mathbb Fe_{2\ov{1}}\oplus \mathbb Fe_{\ov{2} 1},\\
\frg_{\varepsilon_1-\varepsilon_{\ov{1}}}= \mathbb Fe_{1\ov{1}}\oplus \mathbb
Fe_{\ov{2} 2},
 \qquad \frg_{-\varepsilon_1+\varepsilon_{\ov{1}}}
    =\mathbb Fe_{\ov{1} 1}\oplus \mathbb Fe_{2\ov{2}}.
\end{gather*}
Since $\widetilde\frg_\alpha$ projects isomorphically onto
$\frg_\alpha$ for any $\alpha\in\Delta$ (Lemma \ref{uce}),  we
will denote the unique element
 in $\widetilde\frg$ which projects onto $e_{ij}\in\frg$
 (for $i\ne j\in \{ 1,2,\ov{1},\ov{2}\}$)  by $e_{ij}$  also.  Thus
\begin{equation*}
[\widetilde\frg\subo,\widetilde\frg\subo]
 =\espan{e_{12},e_{21},[e_{12},e_{21}],
   e_{\ov{1}\,\ov{2}},e_{\ov{2}\,\ov{1}},
   [e_{\ov{1}\, \ov{2}},e_{\ov{2}\,\ov{1}}]}
\end{equation*}
and
\begin{equation*}
\varepsilon_1([e_{12},e_{21}])=1
 =\varepsilon_{\ov{1}}([e_{\ov{1}\,\ov{2}},e_{\ov{2}\,\ov{1}}]),\quad
\varepsilon_{\ov{1}}([e_{12},e_{21}])=0
 =\varepsilon_1([e_{\ov{1}\,\ov{2}},e_{\ov{2}\,\ov{1}}]).
\end{equation*}
\smallskip

Relative to the system  $\Pi=\{ \alpha=2\varepsilon_1,
\beta=-\varepsilon_1+\varepsilon_{\ov{1}}\}$ of simple roots,
$\varepsilon_1+\varepsilon_{\ov{1}} = \alpha + \beta$ and
$2\varepsilon_{\ov{1}} = \alpha+2\beta$.    Consequently,
$\widetilde\frg$ can be graded  by the `height in $\alpha$',  so
that
\begin{equation*}
\widetilde\frg=\widetilde\frg(-1)\oplus\widetilde\frg(0)\oplus
\widetilde\frg(1),
\end{equation*}
where
\begin{equation*}
\begin{split}
&\widetilde\frg(-1) =\widetilde\frg_{-2\varepsilon_1}\oplus
     \widetilde\frg_{-(\varepsilon_1+\varepsilon_{\ov{1}})}\oplus
     \widetilde\frg_{-2\varepsilon_{\ov{1}}}, \\
&\widetilde\frg(0) =
  \widetilde\frg_{\varepsilon_1-\varepsilon_{\ov{1}}}
  \oplus\widetilde\frg_0\oplus\widetilde\frg_{-
\varepsilon_1+\varepsilon_{\ov{1}}},\\
&\widetilde\frg(1)=\widetilde\frg_{2\varepsilon_1}\oplus
   \widetilde\frg_{\varepsilon_1+\varepsilon_{\ov{1}}}
   \oplus \widetilde\frg_{2\varepsilon_{\ov{1}}}.
\end{split}
\end{equation*}
Set  $h=[e_{12},e_{21}]+[e_{\ov{1}\,\ov{2}},e_{\ov{2}\,\ov{1}}]$,
\ $e=e_{12}+e_{\ov{1}\,\ov{2}}$, and $f=e_{21}+e_{\ov{2}\,\ov{1}}$.
Then $h,e,f$ span a copy of $\frsl_2$ inside
$[\widetilde\frg\subo,\widetilde\frg\subo]$ (the ``diagonal'' subalgebra
of $[\widetilde\frg\subo,\widetilde\frg\subo]=\frsl_2\oplus\frsl_2$).
The gradation spaces above can be described alternatively by
\begin{equation*}
\widetilde\frg(i)=\{ x\in \widetilde\frg \mid [h,x]=2ix\},\quad i=0,\pm 1.
\end{equation*}

Assume now that $L$ is an A$(1,1)$-graded Lie superalgebra  with grading
subsuperalgebra
$\widetilde\frg$, a central cover of $\psl(2,2)$.   Then $L$ also has a
decomposition:
\begin{equation*}
L=L(-1)\oplus L(0)\oplus L(1),
\end{equation*}
where
\begin{equation*}
\begin{split}
&L(-1) =L_{-2\varepsilon_1}\oplus
     L_{-(\varepsilon_1+\varepsilon_{\ov{1}})}\oplus
     L_{-2\varepsilon_{\ov{1}}}=\{x\in L\mid [h,x]=-2x\}, \\
&L(0) =
  L_{\varepsilon_1-\varepsilon_{\ov{1}}}
  \oplus L_0\oplus L_{-\varepsilon_1+\varepsilon_{\ov{1}}}
  =\{x\in L\mid [h,x]=0\},\\
 &L(1)=L_{2\varepsilon_1}\oplus
     L_{\varepsilon_1+\varepsilon_{\ov{1}}}\oplus
     L_{2\varepsilon_{\ov{1}}}=\{x\in L\mid [h,x]=2x\}.
\end{split}
\end{equation*}

Observe that
\begin{equation*}
\begin{split}
L_{\varepsilon_1-\varepsilon_{\ov{1}}}&=
   [[e_{12},e_{21}],L_{\varepsilon_1-\varepsilon_{\ov{1}}}]\\
   &=[e_{12},[e_{21},L_{\varepsilon_1-\varepsilon_{\ov{1}}}]]
    \qquad (\hbox{\rm as} \ \ [e_{12},L_{\varepsilon_1-
\varepsilon_{\ov{1}}}]\subseteq
       L_{3\varepsilon_1-\varepsilon_{\ov{1}}}=0)\\
   &\subseteq [L_{2\varepsilon_1},L_{-(\varepsilon_1+\varepsilon_{\ov{1}}})]
   \subseteq [L(1),L(-1)],
\end{split}
\end{equation*}
and also
\begin{equation*}
\begin{split}
 & L_{-\varepsilon_1+\varepsilon_{\ov{1}}}\subseteq [L(1),L(-1)],\\
 &[L_{\varepsilon_1-\varepsilon_{\ov{1}}},L_{-
\varepsilon_1+\varepsilon_{\ov{1}}}]
 \subseteq [L(0),[L(1),L(-1)]]\subseteq [L(1),L(-1)].
\end{split}
\end{equation*}
Therefore, since
$L_0=\sum_{\gamma\in\Delta}[L_\gamma,L_{-\gamma}]$,  it follows
that
\begin{equation*}
L(0)=[L(1),L(-1)],
\end{equation*}
Consequently,  by the Tits-Kantor-Koecher construction  for Jordan
superalgebras (see, for instance the last paragraph of Section 3,
\cite[Thm. 1.10 and 1.11]{BZ} and the references therein, or
\cite{GN}) $J=L(1)$ is a Jordan superalgebra with multiplication
given by
\begin{equation}\label{Jordan}
x\cdot y =\frac{1}{2}[[x,f],y]
\end{equation}
and unit element $e\, (=e_{12}+e_{\ov{1}\,\ov{2}})$.   Moreover,
$\frj\,:=\widetilde\frg(1)$ is a
Jordan subsuperalgebra of $J$ under this product.  Let
\begin{equation}\label{M11}
e_1:=e_{12},\quad e_2:=e_{\ov{1}\,\ov{2}},\quad x: =e_{1\ov{2}},\quad \hbox{\rm
and}
\quad y:=e_{\ov{1} 2}.
\end{equation}
Then $\{e_1,e_2,x,y\}$ is an homogeneous basis of $\frj$ and,
relative to  the multiplication in \eqref{Jordan}, $\frj$ is
(isomorphic to) the Jordan superalgebra $\mathfrak M: =
\mathfrak{M}_{1,1}(\mathbb F)^+ = \mathfrak M\subo  \oplus
\mathfrak M\subu$, where $\mathfrak M\subo = \big \{
\left(\begin{smallmatrix} a & 0 \\ 0 & d \end{smallmatrix}\right)
\big | \, a,d \in \mathbb F\big \}$, and  $\mathfrak M\subu = \big
\{ \left(\begin{smallmatrix} 0 & b \\ c & 0
\end{smallmatrix}\right) \big | \,$ $ b,c\in \mathbb F\big \}$.  The
product in $\mathfrak M$ is simply the (super)symmetrized matrix
multiplication

\begin{equation}\label{Jproduct}
X\cdot Y = \frac{1}{2}\Big (XY + (-1)^{|X|\,|Y|} YX \Big ). \end{equation}
The
isomorphism can be realized explicitly by matching the basis
elements of $\frj$ with matrix units in $\mathfrak M$ according to

\begin{equation}\label{matrixu}
e_1= \left(\begin{smallmatrix} 1&0\\ 0&0\end{smallmatrix}\right),
\quad e_2=\left(\begin{smallmatrix} 0&0\\
    0&1\end{smallmatrix}\right),
\quad x=\left(\begin{smallmatrix} 0&1\\
    0&0\end{smallmatrix}\right),
\quad y=\left(\begin{smallmatrix} 0&0\\
     1&0\end{smallmatrix}\right).
\end{equation}

\noindent Thus,
\begin{gather} \label{M11prod} e_1^2 = e_1,
\qquad e_2^{2} = e_2, \qquad e_1\cdot e_2= 0 \\
x\cdot y = e_1-e_2= -y\cdot x, \qquad e_1\cdot x = \frac{1}{2} x =
e_2\cdot x, \qquad e_1\cdot y =  \frac{1}{2} y = e_2\cdot y.
\nonumber\end{gather}

Therefore, by the same arguments as in \cite[1.14]{BZ}, $L$ is
centrally isogenous to the centerless Tits-Kantor-Koecher Lie
superalgebra $\TKK(J)$, and the first part of the next (and last)
theorem is established:

\medskip

\begin{theorem}
The A$(1,1)$-graded Lie superalgebras are precisely those Lie
superalgebras which are centrally isogenous to a
Tits-Kantor-Koecher Lie superalgebra $\TKK(J)$  of a Jordan superalgebra $J$
which contains $\mathfrak{M}_{1,1}(\mathbb F)^+$ as a unital subsuperalgebra.
\end{theorem}
\begin{proof}
Only the converse remains to be proven. Let $J$ be a Jordan
superalgebra.  Set  $\frj=\mathfrak{M}_{1,1}(\mathbb F)^+$, and
assume that $\frj$ is a unital subsuperalgebra of $J$.  It has to
be shown that $T=\TKK(J)$ is  A$(1,1)$-graded.  Recall that $T =
T(-1) \oplus T(0) \oplus T(1)$, where $T(1)=J$ and
$T(-1)=\overline J$ is another copy of $J$.   Consider the Jordan
supersystem $\uJ=(T(i),i=\pm 1)$.   Also, let
$\frt=\TKK(\frj)=\frt(-1)\oplus\frt(0)\oplus\frt(1)$, which is
isomorphic to $\psl(2,2)$, and let $\ufrj=(\frt(i), i=\pm 1)$ be
the corresponding Jordan supersystem.   The inclusion
$\frj\subseteq J$ gives a homomorphism of Jordan supersystems so,
by Section 3, there is a unique homomorphism of $3$-graded Lie
superalgebras
\begin{equation*}
\nu : U_{\ufrj}=U_{\ufrj}({-1})\oplus U_{\ufrj}(0) \oplus
U_{\ufrj}(1) \rightarrow T.
\end{equation*}
Moreover, the universal Lie superalgebra $U_{\ufrj}$ is centrally
closed and centrally isogenous to $\frt$, and $\ker\nu\subseteq
U_\ufrj(0)$.  Hence $[\ker\nu,\frt(\pm 1)]=0$, so
$\ker\nu\subseteq Z(U_\ufrj)$. Thus $\nu(U_{\ufrj})\cong
U_{\ufrj}/\ker\nu$ is a central cover of $\frt$ ($\frt$ is
centerless).

Therefore $T$ contains $\widetilde\frg=\nu(U_{\ufrj})$, which is a
central cover of $\frg$, and
$\widetilde\frg=\widetilde\frg(-
1)\oplus\widetilde\frg(0)\oplus\widetilde\frg(1)$
as above, where $\widetilde\frg(1)=\frj$.

For $a \in J=T(1)$, let $\ov{a}$ denote the copy of
$a$ in $\overline J = T(-1)$.  Then for $a,b,c \in J$ we have
\begin{equation}\label{abbc}
[[a,\ov{b}],c]=2\bigl((a\cdot b)\cdot c+a\cdot (b\cdot c)-b\cdot
(a\cdot c)\bigr),
\end{equation}
while
\begin{equation}\label{abbbc}
[[a,\ov{b}], \ov{c}]=2\overline{\bigl(-(a\cdot b)\cdot c+a\cdot
(b\cdot c)-b\cdot (a\cdot c)\bigr)}.
\end{equation}
Take the basis of $\frj=\mathfrak{M}_{1,1}(\mathbb F)^+$ formed by the elements
in (\ref {matrixu}).
Then $e_1$ and $e_2$ are orthogonal idempotents in $J$ with
$e_1+e_2=1$, so $J$ has the Peirce decomposition:
\begin{equation}\label{Peirce}
J=J_0\oplus J_1\oplus J_2,
\end{equation}
where $J_i=\{ v\in J \mid e_1\cdot v=\frac{i}{2}v\}$.

Now, the Lie superalgebra $\frt=\TKK(\frj)$ is isomorphic to
$\psl(2,2)$ via the homomorphism which makes the following assignments:
\begin{align*}
e_1&\mapsto e_{12}& e_2&\mapsto e_{\ov{1}\,\ov{2}}
  & x&\mapsto e_{1\ov{2}} & y&\mapsto e_{\ov{1} 2}\\
\ov{e}_1&\mapsto e_{21}& \ov{e}_2&\mapsto e_{\ov{2}\,\ov{1}}
  & \ov{x}&\mapsto e_{2 \ov{1}} & \ov{y}&\mapsto e_{\ov{2} 1}
\end{align*}
Since $\widetilde\frg(1)=\frt(1)=\frj$ and
$\widetilde\frg(-1)=\frt(-1)=\ov{\frj}$, $\widetilde\frg$ is generated, as
a Lie superalgebra, by the elements $e_1,e_2,x,y,
\ov{e}_1,\ov{e}_2,\ov{x},\ov{y}$,
and $[\widetilde\frg\subo,\widetilde\frg\subo]$ is the
span of $\{ e_1,\ov{e}_1,[e_1,\ov{e}_1],e_2,\ov{e}_2,[e_2,\ov{e}_2]\}$. That is,
\begin{equation*}
[\widetilde\frg\subo,\widetilde\frg\subo]=
 \espan{e_1,\ov{e}_1,[e_1,\ov{e}_1]}\oplus
 \espan{e_2,\ov{e}_2,[e_2,\ov{e}_2]},
\end{equation*}
the direct sum of two copies of $\frsl_2$, with Cartan subalgebra
$\widetilde\frh =\mathbb F[e_1,\ov{e}_1]\oplus \mathbb F[e_2,\ov{e}_2]$.
Therefore, it
suffices to show that relative to the action of $\widetilde\frh$, the Lie
superalgebra $T$
decomposes as
\begin{equation}\label{Tgamma}
T=T_0\oplus\Big(\bigoplus_{\gamma\in \Delta}T_\gamma\Big),
\end{equation}
where
\begin{equation}\label{T0}
T_0=\sum_{\gamma\in\Delta} [T_\gamma,T_{-\gamma}].
\end{equation}
But from \eqref{abbc} and \eqref{abbbc},  it follows that:
\begin{equation}\label{Jis}
\begin{aligned}
J_2&\subseteq T_{2\varepsilon_1},&J_0
  &\subseteq T_{2\varepsilon_{\ov{1}}},
  &J_1&\subseteq T_{\varepsilon_1+\varepsilon_{\ov{1}}},\\
\ov{J}_2&\subseteq T_{-2\varepsilon_1},&\qquad\ov{J}_0
  &\subseteq T_{-2\varepsilon_{\ov{1}}},
  &\qquad\ov{J}_1&\subseteq T_{-(\varepsilon_1+\varepsilon_{\ov{1}})}.
\end{aligned}
\end{equation}
Since $T(0)=[T(1),T(-1)]=[J,\ov{J}]$, the containments in \eqref{Jis}
are in fact equalities  and
\begin{equation*}
T(0)=[J,\overline J]=T_0\oplus T_{2(\varepsilon_1-\varepsilon_{\ov{1}})}
  \oplus T_{-2(\varepsilon_1-\varepsilon_{\ov{1}})}
  \oplus T_{\varepsilon_1-\varepsilon_{\ov{1}}}
  \oplus T_{-\varepsilon_1+\varepsilon_{\ov{1}}}.
\end{equation*}
But $T_{2(\varepsilon_1-\varepsilon_{\ov{1}})}=[J_2,\ov{J}_0]=0$ because
of the usual properties of the Peirce decomposition in
\eqref{Peirce}; namely, \  $J_2\cdot J_0=0=(J_2,J,J_0)$ (where $(\,,\,,\,)$
denotes the associator: $(a,b,c)=(a\cdot b)\cdot c-a\cdot (b\cdot
c)$). Hence \eqref{Tgamma} is obtained, while \eqref{T0} follows
immediately from the condition $T(0)=[T(1),T(-1)]$.
\end{proof}

\bigskip

\begin{remark} \ The results in Section 3 show that, given a
unital Jordan superalgebra $J$ and the associated Jordan
supersystem $\uJ=(J_1,J_{-1})$ (where $J_{\pm 1}=J$), then $U_\uJ$
is the universal central extension of $\TKK(J)$. Therefore, any
A$(1,1)$-graded Lie superalgebra is a central quotient of a
Lie superalgebra $U_\uJ$, for $\uJ$ as above, where $J$ is a
Jordan superalgebra which contains $\mathfrak{M}_{1,1}(\mathbb
F)^+$ as a unital subsuperalgebra.
\end{remark}

\medskip

\begin{remark} \   It is easy to see that the Jordan
superalgebras $\mathfrak{M}_{n,n}(\mathbb F)^+$, $n\geq 2$,
contain $\mathfrak{M}_{1,1}(\mathbb F)^+$, so their
Tits-Kantor-Koecher superalgebras are A$(1,1)$-graded.

The  simple Jordan superalgebra

$$\hbox{\rm JP}_4(\mathbb F) = \left \{ \left(\begin{smallmatrix} a&b\\
   c&a^{\mathfrak t} \end{smallmatrix}\right) \ \big | \ a,b,c\in
\mathfrak{M}_4(\mathbb F),
  \ b = -b^{\mathfrak t} , c = c^{\mathfrak t}\right \} \subset
\mathfrak{M}_{4,4}^+(\mathbb F)$$
 (where $\mathfrak t$ is the transpose) contains $\mathfrak{M}_{1,1}(\mathbb
F)^+$ as a unital subsuperalgebra.
 This containment  may be realized explicitly  in terms of matrix units by the
following
 assignments  (compare the
 related embedding in \cite[Sec. 5]{BE1.5}):

  \begin{eqnarray}\label{JP} e_1 &=& E_{1,1}+E_{2,2} + E_{5,5} + E_{6,6}, \\
   e_2 &=&  E_{3,3}+E_{4,4} + E_{7,7} + E_{8,8},  \nonumber \\
   x &=&  E_{1,7} + E_{2,8} -E_{3,5} - E_{4,6}, \nonumber \\
       y &=&  2\left (E_{7,1}
  +E_{8,2} + E_{5,3}+E_{6,4}\right). \nonumber \end{eqnarray}
Thus, these elements satisfy the multiplication relations in  (\ref{M11prod}).
The Tits-Kantor-Koecher superalgebra of $\hbox{\rm JP}_4(\mathbb F)$ is the
simple  Lie superalgebra P$(3)$, so P$(3)$ is A$(1,1)$-graded.  In fact, much
more
is true.   For an arbitrary unital associative commutative superalgebra
$\Phi$ and an even derivation  $d: \Phi  \rightarrow \Phi$,  the
Cheng Kac Jordan superalgebra $\JCK(\Phi,d)$  contains $\hbox{\rm JP}_4(\mathbb
F)$,
hence also
$\mathfrak{M}_{1,1}(\mathbb F)^+$.    Therefore, the  Cheng-Kac
Lie superalgebra $\CK(\Phi,d)$, which is the Tits-Kantor-Koecher
superalgebra of  $\JCK(\Phi,d)$,  (see \cite{CK,GLS,MZ.Jordan})   is
A$(1,1)$-graded.   As every P$(3)$-graded Lie superalgebra
is centrally isogenous to some $\CK(\Phi,d)$ by  \cite{MZ}, every P$(3)$-graded
superalgebra is A$(1,1)$-graded too.

The simple Jordan superalgebra

$$\hbox{\rm JQ}_4(\mathbb F) =\left \{\left(\begin{smallmatrix} a&b\\
   b&a\end{smallmatrix}\right) \ \big | \ a,b \in \mathfrak{M}_4(\mathbb
F)\right \}
   \subset \mathfrak{M}_{4,4}^+(\mathbb F)$$
   also contains $\mathfrak{M}_{1,1}(\mathbb F)^+$ as a unital subsuperalgebra,
   for example by taking  $e_1$ and $e_2$ as in (\ref{JP})  and
  \begin{eqnarray}\label{JQ}
   x &=&  E_{1,7} + E_{2,8} +E_{5,3} + E_{6,4}, \nonumber \\
       y &=&  2\left (E_{7,1}
  +E_{8,2} + E_{3,5}+E_{4,6}\right), \nonumber \end{eqnarray}
so its Tits-Kantor-Koecher Lie superalgebra, which is Q$(3)$,  is
A$(1,1)$-graded.
\end{remark}

\providecommand{\bysame}{\leavevmode\hbox
to3em{\hrulefill}\thinspace}
\providecommand{\MR}{\relax\ifhmode\unskip\space\fi MR }

\end{document}